\documentclass[12pt]{article}
\usepackage{amsmath}
\usepackage{amssymb}
\usepackage{algorithm}
\usepackage{algorithmic}
\usepackage[dvips]{graphicx}

\def\eps{\varepsilon}
\font\tencmmib=cmmib10 \skewchar\tencmmib '60
\newfam\cmmibfam
\textfont\cmmibfam=\tencmmib

\def\lessim{\ \lower4pt\hbox{$
\buildrel{\displaystyle <}\over\sim$}\ }
\def\gessim{\ \lower4pt\hbox{$\buildrel{\displaystyle >}
\over\sim$}\ }

\def\eps{{\varepsilon}}

\newcommand{\Reals}{\mathbb{R}}

\newtheorem{lemma}{Lemma}
\newtheorem{theorem}{Theorem}

\makeatletter
\@addtoreset{equation}{section}

\makeatother

\def\qed{\hfill$\square$\medskip}

\begin{document}

\title{A note on Talagrand's concentration inequality\\ for
empirical processes}

\author{
Dmitriy Panchenko\thanks{
Department of Mathematics and Statistics, The University of New Mexico,
panchenk@math.unm.edu}
}
\date{January, 2001}

\maketitle
\begin{abstract}
In this paper we revisit Talagrand's proof of concentration
inequality for empirical processes. 
We give a different shorter proof of the main technical lemma
that guarantees the existence of a certain kernel.
Our proof provides an almost optimal value of the constant
involved in the statement of this lemma.
\end{abstract}

\vskip 5mm


\section{Introduction and the proof of main lemma}

This paper was motivated by Section 4 of the
``New concentration inequalities in product spaces''
by Michel Talagrand.
For the most part we will keep the same notations 
with possible minor changes. We slightly weaken the
definition of the distance $m(A,x)$ below compared to \cite{Ta},
but, essentially, this is what is used in the proof
of the concentration inequality for empirical processes. 
Theorem 1 below is Theorem 4.2 in \cite{Ta}
and we assume that the reader is familiar with the proof.
The main technical step, Proposition 4.2 in \cite{Ta}, is proved
differently and constitutes the statement of Lemma 1 below.

Let $n\geq 1$, and $\Omega^n$ be a measurable product space with a product
measure $\mu^n .$ Consider a probability measure
$\nu$ on $\Omega^n$ and fix $x\in \Omega^n.$ 
For $i\leq n,$ if ${\cal C}_i:=\{y\in\Omega^n : y_i\not= x_i\},$ 
we consider the image of the restriction of $\nu$
to ${\cal C}_i$ by the map $y\to y_i,$ and let $d_i$ be its
Radon-Nikodym derivative with respect to $\mu$. In other words, for any function $g\colon \Omega\to\Reals$,
\begin{equation}
\int_{{\cal C}_i}g(y_i)\,d\nu(y)=
\int_{\Omega}g(y_i)d_i(y_i)\,d\mu(y_i).
\label{RND}
\end{equation}
As in \cite{Ta} we  assume that $\Omega$ is finite and each
point is measurable with a positive measure.
Let $m$ be a number of atoms in $\Omega$ and let
$p_1,\ldots,p_m$ be their probabilities.
Consider the function
\[
\psi(x):=\left\{
\begin{array}{cl}
x^2/4, & \mbox{when $x\leq 2,$} \\
x-1, & \mbox{when $x\geq 2,$}
\end{array}
\right.
\]
and, for any set $A\subseteq \Omega^n$, define
$$
m(A,x):=\inf\Bigl\{m(\nu,x) : \nu(A)=1\Bigr\}\quad
\mbox{ where }\quad
m(\nu,x):=\sum_{i\leq n}\int\psi(d_i)d\mu.
$$
\begin{theorem}
Let $L\geq 1.12$ and $P:=\mu^n.$ Then, for any $A\subseteq \Omega^n$,
\begin{equation}
\int\!\exp
\frac{1}{L}m(A,x)
\,dP(x)\leq \frac{1}{P(A)}.
\end{equation}
\end{theorem}

As we mentioned above, the proof is identical to \cite{Ta}
with Proposition 4.2 substituted by the following lemma.

\begin{lemma}
Let $L\geq 1.12$ and $g_1\geq g_2\geq\ldots\geq g_m>0.$ 
There exist 
$\{k_j^i : 1\leq j<i\leq m\}$ 
such that
\begin{equation}
k_j^i\geq 0,\,\,\,\,\,\,\, \sum_{j<i}k_j^i p_j \leq 1
\label{constraints}
\end{equation}
and
\begin{equation}
\sum_{i\leq m}\frac{p_i}{g_i}\exp\Bigl\{
\sum_{j<i}\Bigl(
\log\frac{g_i}{g_j} k_j^i + \frac{1}{L}
\psi(k_j^i)
\Bigr)p_j
\Bigr\}\leq\frac{1}{p_1 g_1+\ldots+p_m g_m}.
\label{kernel}
\end{equation}
\end{lemma}
{\it Remark.} This lemma does not hold for $L\leq 1.07$
(it is easy to construct the counterexample 
for $m=2$),
which means that $L=1.12$ is close to the optimal.

\medskip
\noindent
{\it Proof.} The proof is by induction on the number of atoms $m.$
The statement of lemma is trivial for $m=1.$ 
Note that, in order to show the existence of $\{k_j^i\}$ 
in the statement, one should
try to minimize the left side of (\ref{kernel}) with respect to
$\{k_j^i\}$ under the constraints (\ref{constraints}).
Note also that each term on the left side of (\ref{kernel})
has its own collection of $k_j^i$ for $j<i$ 
and, therefore, minimization can be performed for each term separately.
We assume that $k_j^i$ are chosen in an optimal way minimizing
the left side of (\ref{kernel}) and it will be convenient to take
among all such optimal choices the one maximizing 
$\sum_{j<i}k_j^i p_j$ for all $i\leq m.$ 
To make the induction step we will start by proving 
the following statement, where we assume that $k_j^i$ correspond
to the specific optimal choice indicated above.

\noindent
{\bf Statement.} {\it For all $i\leq m,$ we have
\begin{equation}
\sum_{j<i}k_j^i p_j<1 \Longleftrightarrow 
\log \frac{g_1}{g_i}< \frac{1}{L} \mbox{ and }
\sum_{j<i} 2L\log \frac{g_j}{g_i} p_j < 1.
\label{less}
\end{equation}
In this case 
$k_j^i = 2L \log \frac{g_j}{g_i}.$
}

\noindent
{\it Proof.}  Let us fix $i$ throughout the proof of the statement.
We first assume that the left side of (\ref{less}) holds.
Suppose that 
$\log \frac{g_1}{g_i}\geq \frac{1}{L}.$ 
In this case, since 
$\sup_{x\in\Reals}\psi ' (x)\leq 1,$
one would not increase the left side of (\ref{kernel})
by increasing $k_1^i$ until 
$\sum_{j<i}k_j^i p_j=1$, which contradicts the choice of $k_j^i.$
On the other hand,
$\log \frac{g_1}{g_i}< \frac{1}{L}$ implies that
$\log \frac{g_j}{g_i}< \frac{1}{L}$ and
$k_j^i\leq 2,$ since $\psi(k)=k-1$ for $k\geq 2$
and such choice of $k_j^i$ would only increase the left side of
(\ref{kernel}).
For $k\leq 2,$ $\psi(k)=k^2/4$ and
$$
\mbox{arg}\min\Bigl(
k \log\frac{g_i}{g_j} +\frac{k^2}{4L}
\Bigr)=2L\log\frac{g_j}{g_i}< 2.
$$
(Notice that this proves $\Longleftarrow$, and we continue discussing $\Longrightarrow$.)
If $\sum_{j<i} 2L\log \frac{g_j}{g_i} p_j \geq 1$
then, since we assumed that $\sum_{j<i}k_j^i p_j<1$, the set
$$
{\cal J}:=\{j : 
k_j^i < 2L \log\frac{g_j}{g_i}
\}\not = \emptyset 
$$
is not empty and increasing $k_j^i$ for $j\in {\cal J}$ would decrease the left side
of (\ref{kernel}) --  a contradition. This completes the prove of the statement.
\qed

We now go back to the proof of Lemma 1. Equation (\ref{less}) implies that if
$ \sum_{j<i}k_j^i p_j<1$ then
$\sum_{j<\ell}k_j^\ell p_j<1$ for $\ell\leq i.$
Therefore, the equality
$ \sum_{j<m-1}k_j^{m-1} p_j=1$ 
would imply
$\sum_{j<m}k_j^m p_j=1.$
Let us first consider the case when  
$ \sum_{j<m-1}k_j^{m-1} p_j=1.$
(This step is meaningless for $m=2$ and
should simply be skipped.)
We will now show that 
$k_j^m = k_j^{m-1},$ $j<m-1$ and $k_{m-1}^m = 0.$
Indeed,
\begin{eqnarray}
&&
\inf_{\sum_{j<m}k_j p_j =1}
\sum_{j<m}\Bigl(
\log\frac{g_m}{g_j}k_j+\frac{1}{L}\psi(k_j)
\Bigr)p_j =
\nonumber
\\
&&
=
\log\frac{g_m}{g_{m-1}}+
\inf_{\sum_{j<m}k_j p_j =1}
\Bigl(
\sum_{j<m-1}\Bigl(
\log\frac{g_{m-1}}{g_j}k_j+\frac{1}{L}\psi(k_j)
\Bigr)p_j 
+\frac{1}{L}\psi(k_{m-1})p_{m-1}
\Bigr).
\label{one}
\end{eqnarray}
Since $\sum_{j<m-1}k_j^{m-1} p_j=1$,
it is advantageous to set $k_{m-1}^m=0$
and $k_j^m = k_j^{m-1}$ for $j<m-1 .$
In this case 
$$
\frac{p_m}{g_m}\exp\Bigl\{
\sum_{j<m}\Bigl(
\log\frac{g_m}{g_j} k_j^m + \frac{1}{L}
\psi(k_j^m)
\Bigr)p_j
\Bigr\}=
\frac{p_m}{g_{m-1}}\exp\Bigl\{
\sum_{j<m-1}\Bigl(
\log\frac{g_{m-1}}{g_j} k_j^{m-1} + \frac{1}{L}
\psi(k_j^{m-1})
\Bigr)p_j
\Bigr\}.
$$
By induction assumption, (\ref{kernel})
holds for the sets $(g_1,\ldots,g_{m-1})$ and 
$(p_1,\ldots,p_{m-1}+p_m).$ 
Since $p_{m-1}g_{m-1}+p_m g_m\leq (p_{m-1}+p_m)g_{m-1},$ this implies that
it holds for $(g_1,\ldots,g_{m})$ and 
$(p_1,\ldots,p_m)$, in this case. 

Now we will assume that 
$\sum_{j<m-1}k_j^{m-1} p_j<1$
or, equivalently, 
$\log \frac{g_1}{g_{m-1}}< \frac{1}{L}$ and 
$\sum_{j<m-1} 2L\log \frac{g_j}{g_{m-1}} p_j < 1.$
It is obvious that in this case there exists 
$g_0<g_{m-1}$ such that for 
$g_m\in (g_0, g_{m-1}]$
both $\log \frac{g_1}{g_{m}}< \frac{1}{L}$ and 
$\sum_{j<m} 2L\log \frac{g_j}{g_{m}} p_j < 1$
hold and, therefore,
$\sum_{j<m}k_j^{m} p_j<1.$
We assume that $g_0$ is the smallest number with such
properties.
Let us show that, for fixed $g_1,\ldots,g_{m-1}$, the case of
$g_m < g_0$ can be transformed to $g_m=g_0.$ 
Indeed, take $g_m < g_0$ so that, by (\ref{less}), $\sum_{j<m}k_j^{m} p_j=1.$
In this case (\ref{one}) still holds and implies that
$k_j^m$ do not depend on $g_m$ for $g_m<g_0$ and
$$
\frac{p_m}{g_m}\exp\Bigl\{
\sum_{j<m}\Bigl(
\log\frac{g_m}{g_j} k_j^m + \frac{1}{L}
\psi(k_j^m)
\Bigr)p_j
\Bigr\}=
\frac{p_m}{g_{m-1}}\exp\Bigl\{
\sum_{j<m}\Bigl(
\log\frac{g_{m-1}}{g_j} k_j^{m} + \frac{1}{L}
\psi(k_j^{m})
\Bigr)p_j
\Bigr\}.
$$
This means that for $g_m<g_0$ the left side of
the inequality (\ref{kernel}) does not depend on $g_m.$
Since $(p_1 g_1+\ldots+p_m g_m)^{-1}$ decreases in $g_m$,
it is enough to prove the inequality for $g_m=g_0.$

Hence, we can consider $g_m\in [g_0,g_{m-1}]$ and assume that
$\log \frac{g_1}{g_{m}}\leq \frac{1}{L},$
$\sum_{j<m} 2L\log \frac{g_j}{g_{m}} p_j \leq 1$
and
$k_j^i= 2L\log \frac{g_j}{g_{i}}.$
Note that (\ref{kernel}) can be rewritten as
\begin{equation}
\sum_{i\leq m}\frac{p_i}{g_i}\exp\Bigl\{
-L\sum_{j<i}
\Bigl(\log\frac{g_j}{g_i}\Bigr)^2 p_j
\Bigr\}\leq\frac{1}{p_1 g_1+\ldots+p_m g_m}.
\label{kernel2}
\end{equation}
By induction hypothesis, (\ref{kernel2}) holds for $g_m=g_{m-1}.$
To prove it for $g_m< g_{m-1}$ we will compare the derivatives 
of both sides of (\ref{kernel2}) with respect to $g_m.$
It is enough to have
\begin{eqnarray*}
\frac{p_m}{g_m}\exp\Bigl\{
-L\sum_{j<m}
\Bigl(\log\frac{g_m}{g_j}\Bigr)^2 p_j
\Bigr\}\Bigl(
-\frac{1}{g_m}-2L\sum_{j<m}
\log\frac{g_m}{g_j} p_j \frac{1}{g_m}
\Bigr)
\geq
-\frac{p_m}{(p_1 g_1+\ldots+p_m g_m)^2}
\end{eqnarray*}
or, equivalently,
\begin{eqnarray*}
\exp\Bigl\{
-L\sum_{j<m}
\Bigl(\log\frac{g_m}{g_j}\Bigr)^2 p_j
\Bigr\}\Bigl(
1-2L\sum_{j<m}
\log\frac{g_j}{g_m} p_j
\Bigr)
\leq
\Bigl(\frac{g_m}{p_1 g_1+\ldots+p_m g_m}
\Bigr)^2 .
\end{eqnarray*}
Since $1-x\leq e^{-x}$ for $x\geq 0$, it is enough to show
\begin{eqnarray*}
\exp\Bigl\{
-L\sum_{j<m} p_j
\Bigl(
\Bigl(\log\frac{g_j}{g_m}\Bigr)^2 
+2 \log\frac{g_j}{g_m}
\Bigr)
\Bigr\}
\leq
\Bigl(\frac{g_m}{p_1 g_1+\ldots+p_m g_m}
\Bigr)^2 .
\end{eqnarray*}
One can check that $(\log x)^2 + 2\log x$ is concave
for $x\geq 1.$ If we express 
$g_j=\lambda_j g_1 + (1-\lambda_j)g_m$ for $j=1,\ldots,m-1$
then 
$$
\sum_{j<m} p_j
\Bigl(
\Bigl(\log\frac{g_j}{g_m}\Bigr)^2 
+2 \log\frac{g_j}{g_m}
\Bigr)\geq
\Bigl(\sum_{j<m} p_j\lambda_j\Bigr)
\Bigl(
\Bigl(\log\frac{g_1}{g_m}\Bigr)^2 
+2 \log\frac{g_1}{g_m}
\Bigr),
$$
$$
p_1 g_1 +\ldots+ p_m g_m =
\Bigl(\sum_{j<m} p_j\lambda_j\Bigr)g_1 +
\Bigl(p_m + \sum_{j<m} (1-\lambda_j)p_j\Bigr)g_m.
$$
If we denote $p=\sum_{j<m} p_j\lambda_j$ and 
$t=\log\frac{g_1}{g_m}$, it is enough to prove
\begin{equation}
\exp\Bigl\{
-Lp(t^2 +2t)
\Bigr\}
\leq
\Bigl(\frac{1}{p e^t+1-p}
\Bigr)^2,\,\,\,
0\leq p\leq 1,\,
0\leq t\leq \frac{1}{L}.
\end{equation}
Equivalently,
$$
\varphi(p,t)=
(pe^t +1 -p)
\exp\Bigl\{-\frac{L}{2}p(t^2 +2t)\Bigr\}\leq 1,\,\,\,
0\leq p\leq 1,\,
0\leq t\leq \frac{1}{L}.
$$
We have
$$
\varphi_t^{\prime}(p,t)=\varphi(p,t)\Bigl(
\frac{pe^t}{pe^t +1-p}-Lp(t+1)
\Bigr).
$$
Since $\varphi(p,0)=1$  for all $p>0$, we need 
$\varphi_t^{\prime}(p,0)=p(1-L)\leq 0,$
or $L\geq 1,$ which holds if $L\geq 1.12.$
It is easy to see that 
$\varphi_t^{\prime}(p,t)=0$ at most at one point $t.$
In combination with $\varphi_t^{\prime}(p,0)\leq 0$
it implies that for a fixed $p$ maximum of $\varphi(p,t)$
is attained at $t=0$ or $t=1/L.$ Therefore, we have to
show $\varphi(p,1/L)\leq 1,\,0\leq p\leq1.$
We have,
$$
\varphi_p^{\prime}(p,\frac{1}{L})=
\varphi(p,\frac{1}{L})\Bigl(
\frac{e^{\frac{1}{L}}-1}{pe^{\frac{1}{L}}+1-p}
-\frac{L}{2}\Bigl(
\frac{1}{L^2}+2\frac{1}{L}
\Bigr)
\Bigr).
$$
Since $\varphi(0,\frac{1}{L})=1$ we should have
$\varphi_p^{\prime}(0,\frac{1}{L})\leq 0$
which would also imply
$\varphi_p^{\prime}(p,\frac{1}{L})\leq 0,\,p>0.$
One can check that
$$
\varphi_p^{\prime}(0,\frac{1}{L})=
e^{\frac{1}{L}}-1-\frac{1}{2}\Bigl(
\frac{1}{L}+2 \Bigr)
< 0
$$
for $L\geq 1.12.$
This finishes the proof.
\qed

\section{One concentration inequality for empirical processes}

Given Theorem 1
one can proceed as in \cite{Ta} to obtain
the classical form of concentration inequality
for the empirical process around its mean. 

We will now show that in one special case, 
which allows certain simplifications,
the technique of Talagrand yields a
rather sharp concentration result with explicit constants.
Consider a countable class of measurable functions 
${\cal F}=\{f:\Omega\to [0,1]\}.$
Consider the following function on $\Omega^n$,
$$
Z(x):=\sup_{f\in{\cal F}}\sum_{i\leq n}(\mu f -f(x_i)),\,\,\,
$$
where $\mu f:=\int f d\mu.$
In applications (see e.g. \cite{KP1}, \cite{KP2}, \cite{KPL}), especially in the case when the empirical process is defined over a family of sets,
the uniform variance is often simply bounded by the uniform second moment,
\begin{equation}
n\sup_{f\in {\cal F}}\mathrm{Var}(f)\leq
\sigma^2:=n\sup_{f\in {\cal F}}\mu f^2,
\label{sigma}
\end{equation}
for which one has some apriori bound.
Talagrand's technique gives in this case
the following concentration inequalities.

\begin{theorem}
Let $L=1.12$, $P=\mu^n$ and $M$ be a median of $Z.$
Then
$$
P(Z\geq M + 2\max(Lu,\sigma\sqrt{L u}))\leq
2e^{-u},
$$
$$
P(Z\leq M - 2\max(Lu,\sigma\sqrt{L u}))\leq
2e^{-u}.
$$
\end{theorem}
{\it Proof.}
Without loss of generality we assume that $\cal F$ is finite.
Given $a\in \Reals$, let us
consider the set $A=\{Z(x)\leq a\}.$
For a fixed $x$, let $f\in{\cal F}$ be such that
\begin{equation}
Z(x)= \sum_{i\leq n}(\mu f - f(x_i)).
\label{sup}
\end{equation}
Then, for any probability measure $\nu$ such that
$\nu(A)=1$, we can write (recall (\ref{RND}))
\begin{eqnarray*}
&&
Z(x)-a\leq
\int \Bigl(
\sum_{i\leq n}(\mu f - f(x_i))-
\sum_{i\leq n}(\mu f - f(y_i))
\Bigr)d\nu(y)
\\
&&
=
\sum_{i\leq n}\int 
(f(y_i)-f(x_i))d_i(y_i)d\mu(y_i)
\leq
\sum_{i\leq n}\int 
f(y_i)d_i(y_i)d\mu(y_i).
\end{eqnarray*}
As is easily checked, $uv\leq u^2+\psi(v)$ for $v\geq 0$ and
$0\leq u\leq 1.$ Therefore, for any $\delta\geq 1$,
\begin{eqnarray*}
Z(x)-a\leq
\delta \sum_{i\leq n}\int 
\frac{f(y_i)}{\delta}d_i(y_i)d\mu(y_i)\leq
\frac{\sigma^2}{\delta} + 
\delta \sum_{i\leq n}\int \psi(d_i)d\mu. 
\end{eqnarray*}
Taking the infimum over $\nu$ we get
$$
Z(x)\leq 
a + \frac{\sigma^2}{\delta} +
\delta m(A,x). 
$$
Theorem 1 then implies that, for $L=1.12$, with probability at least
$1-\frac{1}{P(Z\leq a)}e^{-u}$,
$$
Z(x)\leq a + 2\max(Lu,\sigma\sqrt{Lu}).
$$
Applying this to $a=M$ and to
$a=M- 2\max(Lu,\sigma\sqrt{Lu})$ gives
the result.
\qed

{\it Remark.} It is interesting to note that the bounds
of Theorem 2 seem to avoid a `singular' behaviour 
of the general bounds expressed in terms of the weak variance 
$n\sup_{f\in {\cal F}}\mbox{Var}f$ (see \cite{Ma}),
when the linear dependance of the term
$(1+\eps)M$ on $\eps$
requires the factor of the order $\eps^{-1}$
in the last term of the bound 
$\eps^{-1}u.$ Under the assumptions of Theorem 2,
one can also avoid this `singularity' using a recent result of 
Emmanuel Rio \cite{Rio}, that provides rather sharp constants too,
and the concentration is around mean instead of median.

{\bf Acknowledgments.} We want to thank Michel Talagrand
for pointing out a recent result of Emmanuel Rio.

\end{document}